\documentclass[twocolumn]{article}

\usepackage{microtype}
\usepackage{graphicx}
\usepackage{subcaption}
\usepackage{booktabs}
\usepackage{stmaryrd}
\usepackage{amsmath}
\usepackage{amssymb}
\usepackage{mathtools}
\usepackage{amsthm}
\usepackage{hyperref}
\usepackage[capitalize,noabbrev]{cleveref}
\usepackage{xcolor}
\usepackage{listings}
\usepackage{natbib}
\usepackage{algorithm}
\usepackage{algorithmic}

\theoremstyle{plain}
\newtheorem{theorem}{Theorem}[section]

\theoremstyle{definition}
\newtheorem{definition}[theorem]{Definition}

\theoremstyle{remark}

\definecolor{codegreen}{rgb}{0,0.6,0}
\definecolor{codegray}{rgb}{0.5,0.5,0.5}
\definecolor{codepurple}{rgb}{0.58,0,0.82}
\definecolor{backcolour}{rgb}{0.95,0.95,0.92}
\definecolor{bluekeywords}{rgb}{0,0,1}

\lstdefinelanguage{lean}{morekeywords={def, theorem, lemma, example, definition, structure, instance,
		inductive, match, with, where, let, in, if, then, else,
		forall, fun, $\lambda$, have, show, from, by, begin, end,
		import, open, namespace, section, variable, variables,
		noncomputable, @[simp], @[ext], @[class], @[instance],
		deriving, extends, of, sorry, axiom, constant,
		protected, private, mutual, abstract},
	sensitive=true,
	morecomment=[l]{--},
	morecomment=[s]{/-}{-/},
	morestring=[b]",
	morestring=[b]',
	morestring=[b]"""}  
\lstdefinestyle{leanstyle}{backgroundcolor=\color{backcolour},
	frame=single,
	rulecolor=\color{black},
	commentstyle=\color{codegreen},
	keywordstyle=\color{bluekeywords},
	numberstyle= \tiny\color{codegray},
	stringstyle=\color{codepurple},
	basicstyle=\ttfamily\footnotesize,
	breakatwhitespace=false,
	breaklines=true,
	captionpos=b,
	keepspaces=true,
	numbers=left,
	numbersep=5pt,
	showspaces=false,
	showstringspaces=false,
	showtabs=false,
	tabsize=2,
	alsoletter={_},
	columns=fullflexible,
	mathescape=true, 
}
\title{Formalization of Line Search Methods by Lean}

\author{Yiyang Zhang \\
	School of Data Science, \\
	The Chinese University of Hong Kong, Shenzhen \\
	\texttt{yiyangzhang1@link.cuhk.edu.cn}
	\and
	Kenneth SHUM \\
	School of Science and Engineering, \\
	The Chinese University of Hong Kong, Shenzhen \\
	\texttt{wkshum@cuhk.edu.cn}}

\begin{document}

\maketitle

\begin{abstract}
This paper presents a formalization of line search methods in the Lean 4 theorem prover. Our goal is to advance  machine verification of nonlinear optimization theory by translating standard textbook definitions and convergence arguments into rigorous Lean code. We formalize fundamental notions related to gradient descent and descent directions, adaptive step-size selection via backtracking line search, and several classical line search criteria, including the Armijo, Goldstein, and Wolfe conditions, as well as nonmonotone variants. We further formalize a key convergence result, namely the Zoutendijk theorem, which plays a central role in the global convergence analysis of gradient-based iterative methods. By providing machine-checkable definitions and proofs for line search theory, this work complements existing formalizations of first-order optimization methods and establishes a foundation for the verified development of more advanced algorithms in nonlinear programming.
\end{abstract}

\section{Introduction}

Lean 4~\cite{demoura2021lean4} is both a programming language and an interactive theorem prover designed for formalizing mathematics and verifying programs. It helps bridge the gap between mathematical proofs and executable code by representing proof developments as machine-checkable and shareable programs. In a typical Lean workflow, users write proof scripts in an editor while simultaneously inspecting the current proof states and goals. Lean's elaborator and automation can simplify complex derivations, and users can draw on large searchable theorem libraries such as Mathlib~\cite{mathlib2020}, the largest and most widely used library in the Lean community. In this paper we consider the verification of the correctness of line search method in optimization. For the formalization of other types of optimization algorithms in Lean, we refer the readers to the work of the Beijing International Center for Mathematical Research~\cite{bicmr2026} at Peking University.

The authors of~\cite{li2024formalization} formalize
first-order algorithms for convex optimization, including gradient descent, subgradient descent, proximal gradient, Nesterov's accelerated method, with fully verified convergence rate guarantees. The framework has been further extended to block-structured optimization~\cite{li2025block}, covering block coordinate descent (BCD), alternating direction method of multipliers (ADMM), subdifferentials, and the Kurdyka–Łojasiewicz (KL) property for nonconvex convergence analysis. Moreover, first-order optimality conditions for smooth constrained optimization have also been formalized~\cite{li2025kkt}, including Karush–Kuhn–Tucker (KKT) conditions, linear independence constraint qualification (LICQ), the Farkas lemma, and weak duality. These efforts establish a rigorous, machine-verified foundation for optimization theory in Lean and align closely with the formalization of line search methods presented in this work.

Early optimization methods adopted exact line search techniques, such as interval elimination and quadratic interpolation, to locate the exact one-dimensional minimum at each iteration \cite{kiefer1953sequential,curry1944method}. However, exact line search is numerically unstable in practical optimization. With the development of quasi-Newton methods in the 1960s \cite{davidon1959variable,fletcher1963rapidly}, traditional exact line search faced a bottleneck. To address this problem, Armijo established the basic theory of modern inexact line search in 1966 by proposing the sufficient decrease condition and a backtracking step-size strategy under the Lipschitz gradient assumption, which greatly improved the efficiency and practicability of iterative optimization algorithms \cite{armijo1966}.

In this paper, we formalize fundamental definitions and lemmas for optimization line search theory in Lean~4. Our development covers the adaptive step-size selection rule of backtracking line search, classical line search stopping criteria including the Armijo condition~\cite{armijo1966}, Goldstein condition~\cite{goldstein1965}, and Wolfe conditions~\cite{wolfe1969,wolfe1971}, as well as two non-monotonic conditions \cite{grippo1986nonmonotone,zhang2004nonmonotone}. Lastly, we formulate and prove the pivotal Zoutendijk theorem for iterative optimization convergence analysis \cite{zoutendijk1960,yang2012}. This machine-verified formalization supplements existing formal optimization frameworks, laying a solid foundation for the formal verification of gradient-based iterative optimization algorithms.

\paragraph{Paper Organization.}
The rest of the paper is organized as follows. Section~\ref{sec:math_foundation} briefly reviews the mathematical foundations, which are respectively vector space modelling, inner product and norm, and basic assumptions about gradient. Section~\ref{sec:descent_line_search} introduces the descent and line search framework, including the Armijo, Goldstein, Wolfe, Grippo, and Zhang-Hager conditions. Lastly, Section~\ref{sec:convergence} presents the deviation and the Lean formalization of the Zoutendijk theorem and its corollary.

\section{Mathematical Foundations}
\label{sec:math_foundation}

\subsection{Vector Space, inner product, and norm}

We model the space as a finite-dimensional real vector space.
Specifically, for a fixed dimension $n \in \mathbb{N}$, we represent vectors as functions from the finite index set to the real numbers:
\[
x \in \mathbb{R}^n.
\]
In Lean, we represent it as a function from $\texttt{Fin n}$ to $\mathbb{R}$, where $\texttt{Fin n}$ is the type of consisting of $n$ elements.
\begin{lstlisting}[
	basicstyle=\ttfamily\scriptsize,
	backgroundcolor=\color{gray!5},
	frame=single,
	mathescape=true,
	escapeinside={(*@}{@*)}]
(*@\textcolor{red}{abbrev}@*) (*@\textcolor{blue}{Vec}@*) (n : (*@$\mathbb{N}$@*)) := Fin n (*@$\rightarrow$@*) (*@$\mathbb{R}$@*)
\end{lstlisting}

We consider an objective function
\[
f : \mathbb{R}^n \to \mathbb{R},
\]
which will be the target of our optimization procedure, and corresponding Lean declaration is as follows:

\begin{lstlisting}[
	basicstyle=\ttfamily\scriptsize,
	backgroundcolor=\color{gray!5},
	frame=single,
	mathescape=true,
	escapeinside={(*@}{@*)}]
(*@\textcolor{red}{variable}@*) (f : Vec n (*@$\rightarrow$@*) (*@$\mathbb{R}$@*))
\end{lstlisting}


In order to equip the vector space with a notion of length, we define the Euclidean inner product
\[
\langle x,y\rangle
:=
\sum_{i=1}^{n} x_i y_i.
\]
This inner product is positive semi-definite,
\[
\langle x, x \rangle \ge 0 \quad \text{for all } x,
\]
which we do not postulate but \emph{prove}: it follows from the term-wise non-negativity $x_i^2 \ge 0$ together with \texttt{Finset.sum\_nonneg}.

The induced \emph{Euclidean} norm is defined as
\[
\|x\|
:=
\sqrt{\langle x,x\rangle},
\]
for $x$ in $\texttt{Vec n}$.
We emphasize that this is a custom Euclidean norm \texttt{vnorm}, and \emph{not} the ambient norm $\|\cdot\|$ that Mathlib places on $\texttt{Fin n} \to \mathbb{R}$, which is the \emph{sup} norm.
For the sup norm, both the identity $\|x\|^2 = \langle x,x\rangle$ and the Cauchy–Schwarz inequality \emph{fail} for $n \ge 2$, which is precisely why an earlier version had to postulate them as axioms.
With the genuine Euclidean norm \texttt{vnorm}, both facts become \emph{theorems} rather than assumptions:
\[
\|x\|^2 = \langle x,x\rangle
\]
follows from \texttt{Real.sq\_sqrt}, and the Cauchy–Schwarz inequality
\[
\langle x,y\rangle
\le
\|x\| \|y\|
\]
holds for all $x$ and $y$ in $\texttt{Vec n}$ by Mathlib's finite Cauchy–Schwarz inequality \texttt{Finset.sum\_mul\_sq\_le\_sq\_mul\_sq}.
We record them as the following definitions and lemmas in Lean:
\begin{lstlisting}[
	basicstyle=\ttfamily\scriptsize,
	backgroundcolor=\color{gray!5},
	frame=single,
	mathescape=true,
	escapeinside={(*@}{@*)}
	]
(*@\textcolor{red}{def}@*) (*@\textcolor{blue}{innerpro}@*) (x y : Vec n) : (*@$\mathbb{R}$@*) :=
	(*@$\sum$@*) i : Fin n, x i * y i
	
(*@\textcolor{red}{noncomputable def}@*) (*@\textcolor{blue}{vnorm}@*) (x : Vec n) : (*@$\mathbb{R}$@*) :=
	Real.sqrt (innerpro x x)
	
(*@\textcolor{red}{lemma}@*) (*@\textcolor{blue}{innerproSelfNonneg}@*) (x : Vec n) :
	innerpro x x (*@$\ge$@*) (*@\textcolor{green!50!black}{0}@*)
	
(*@\textcolor{red}{lemma}@*) (*@\textcolor{blue}{vnorm\_sq}@*) (x : Vec n) :
	(vnorm x)^(*@\textcolor{green!50!black}{2}@*) = innerpro x x
	
(*@\textcolor{red}{lemma}@*) (*@\textcolor{blue}{cauchySchwarz}@*) (x y : Vec n) :
	innerpro x y (*@$\le$@*) vnorm x * vnorm y
\end{lstlisting}

\subsection{Gradient Lipschitz Continuity}
We assume the objective function admits a gradient mapping:
\[
\nabla f : \mathbb{R}^n \to \mathbb{R}^n
\]
which is Lipschitz continuous.
\begin{definition}[Lipschitz Continuity]
	{\it A differentiable function $f$ is said to be gradient Lipschitz continuous if there exists a fixed constant $L>0$ such that}
	\[
	\|\nabla f(x)-\nabla f(y)\|
	\le
	L\|x-y\|
	\]
	{\it for all $x$ and $y$, where $\|\cdot\|$ denotes the Euclidean norm \texttt{vnorm}.}
\end{definition}
The assumption of Lipschitz continuity is standard in nonlinear optimization, ensuring that the function is smooth and enabling descent-based methods to behave predictably.
We emphasize that the constant $L$ is a \emph{fixed} parameter of the predicate, and is \emph{not} universally quantified.
An earlier formulation wrote ``$\forall L>0,\dots$'', which asserts the bound for \emph{every} positive $L$; letting $L\to 0$, that degenerate statement would force $\nabla f$ to be constant, so it is not the genuine Lipschitz condition.
Moreover, since this property cannot be proved for an arbitrary pair $(f,\nabla f)$, it is no longer postulated as an axiom: it is carried as an \emph{explicit hypothesis} (a predicate \texttt{GradLipschitz}) of the theorems that need it.
The Lean code is articulated as follows:
\begin{lstlisting}[
	basicstyle=\ttfamily\scriptsize,
	backgroundcolor=\color{gray!5},
	frame=single,
	mathescape=true,
	escapeinside={(*@}{@*)}
	]
(*@\textcolor{red}{variable}@*) (grad : Vec n (*@$\rightarrow$@*) Vec n)
	
(*@\textcolor{red}{def}@*) (*@\textcolor{blue}{GradLipschitz}@*)
	(g : Vec n (*@$\rightarrow$@*) Vec n) (L : (*@$\mathbb{R}$@*)) : Prop :=
	(*@$\forall$@*) x y : Vec n,
	vnorm (g x - g y) (*@$\le$@*) L * vnorm ((*@\textcolor{red}{fun}@*) i => x i - y i)
\end{lstlisting}
To connect the gradient with descent behavior, we assume a first-order approximation property.
\begin{definition}[1st-order approximation property]
	{\it Let $c_1\in(0,1)$ be an Armijo parameter. For any point $x$ and any descent direction $d$, both with type} $\texttt{Vec n}$ {\it (that is, $\langle \nabla f(x),d\rangle<0$), there exists a sufficiently small step size $\varepsilon>0$ such that, for every $\alpha\in(0,\varepsilon)$,}
	\begin{equation}
		f(x+\alpha d)
		\le
		f(x)+c_1\,\alpha\, \langle \nabla f(x),d\rangle.
		\label{eq:first_order_approximation}
	\end{equation}
\end{definition}
This property is the genuine first-order consequence of differentiability.
By the Taylor expansion $f(x+\alpha d)=f(x)+\alpha\langle\nabla f(x),d\rangle+o(\alpha)$, and since $d$ is a descent direction the quantity $(1-c_1)\,\alpha\,\langle\nabla f(x),d\rangle$ is strictly negative, so the Armijo inequality~\eqref{eq:first_order_approximation} holds for all sufficiently small $\alpha>0$.
We stress that the earlier formulation used the coefficient $1$ (with a strict inequality),
\[
f(x+\alpha d) < f(x)+\alpha\,\langle\nabla f(x),d\rangle,
\]
which would require the remainder to satisfy $o(\alpha)<0$; this is in general \emph{false} for smooth (convex) functions, making that statement \emph{unsatisfiable}. Introducing the Armijo parameter $c_1\in(0,1)$ together with the descent-direction premise yields the genuinely satisfiable condition, which we use in the proof of the descent direction.
\begin{lstlisting}[
	basicstyle=\ttfamily\scriptsize,
	backgroundcolor=\color{gray!5},
	frame=single,
	mathescape=true,
	escapeinside={(*@}{@*)}
	]
(*@\textcolor{red}{def}@*) (*@\textcolor{blue}{FirstOrderApprox}@*) 
  (f : Vec n (*@$\rightarrow$@*) (*@$\mathbb{R}$@*)) (g : Vec n (*@$\rightarrow$@*) Vec n) : Prop :=
  (*@$\forall$@*) (x d : Vec n) (c1 : (*@$\mathbb{R}$@*)), (*@\textcolor{green!50!black}{0}@*) < c1 (*@$\rightarrow$@*) c1 < (*@\textcolor{green!50!black}{1}@*) (*@$\rightarrow$@*)
    innerpro (g x) d < (*@\textcolor{green!50!black}{0}@*) (*@$\rightarrow$@*)
    (*@$\exists$@*) (*@$\varepsilon$@*) > (*@\textcolor{green!50!black}{0}@*),
    (*@$\forall$@*) (*@$\alpha$@*) (*@$\in$@*) Set.Ioo ((*@\textcolor{green!50!black}{0}@*) : (*@$\mathbb{R}$@*)) (*@$\varepsilon$@*),
       f ((*@\textcolor{red}{fun}@*) i => x i + (*@$\alpha$@*) * d i)
         (*@$\le$@*) f x + c1 * (*@$\alpha$@*) * innerpro (g x) d
\end{lstlisting}


\section{Descent and Line Search}
\label{sec:descent_line_search}

\subsection{Descent Direction}

\begin{definition}[Descent Direction]
	{\it A direction \(d\) is called a descent direction at point \(x\) if}
	\[
	\langle \nabla f(x),d\rangle < 0.
	\]
\end{definition}

This condition ensures that moving along the direction d could locally decrease the objective function. We can use Lean to formalize it:

\begin{lstlisting}[
	basicstyle=\ttfamily\scriptsize,
	backgroundcolor=\color{gray!5},
	frame=single,
	mathescape=true,
	escapeinside={(*@}{@*)}
	]
(*@\textcolor{red}{def}@*) (*@\textcolor{blue}{isDescentDirection}@*) (x d : Vec n) : (*@\textcolor{red}{Prop}@*) :=
  innerpro (grad x) d < (*@\textcolor{green!50!black}{0}@*)
\end{lstlisting}
Hence, a descent direction guarantees the existence of a step size that makes the function value smaller:
\[
\exists \alpha > 0,\ f(x + \alpha d) < f(x),
\]
which can be defined in Lean as:
\begin{lstlisting}[
	basicstyle=\ttfamily\scriptsize,
	backgroundcolor=\color{gray!5},
	frame=single,
	mathescape=true,
	escapeinside={(*@}{@*)}
	]
(*@\textcolor{red}{lemma}@*) (*@\textcolor{blue}{descentImpliesDecrease}@*)
  (x d : Vec n)
  (hdesc : isDescentDirection grad x d) 
  (hfo : FirstOrderApprox f grad):
  (*@$\exists$@*) $\alpha$ > (*@\textcolor{green!50!black}{0}@*), f ((*@\textcolor{red}{fun}@*) i => x i + $\alpha$ * d i) < f x
\end{lstlisting}

Using \texttt{firstOrderApprox}, we know that for sufficiently small positive step sizes $\alpha$ such that \eqref{eq:first_order_approximation} holds.
We then choose a specific small step size $\alpha=\varepsilon/2$. Since $\alpha>0$ and the inner product is negative, we have
$$\alpha\langle \nabla f(x), d\rangle <0,$$
which implies that the right-hand side of \eqref{eq:first_order_approximation} is strictly smaller than $f(x)$. Therefore,
$f(x+\alpha d)<f(x)$. This shows that there exists a positive step size that yields a strict decrease in the function value.

\subsection{Armijo Line Search}

The Armijo condition only enforces sufficient decrease in the objective function, is easy to implement, and is often paired with backtracking line search, but it tends to produce overly small step sizes. It is suitable for basic gradient descent, and it is defined as follows.

\begin{definition}[Armijo Condition]
	{\it Assume \(d^k\) is a descent direction of \(x^k\). A step size \(\alpha > 0\) satisfies the Armijo condition if:}
	\[
	f(x^k + \alpha d^k) \le f(x^k) + c_1 \alpha \left\langle \nabla f(x^k), d^k \right\rangle, \quad c_1 \in (0,1).
	\]
\end{definition}

\begin{lstlisting}[
	basicstyle=\ttfamily\scriptsize,
	backgroundcolor=\color{gray!5},
	frame=single,
	mathescape=true,
	escapeinside={(*@}{@*)}
	]
(*@\textcolor{red}{def}@*) (*@\textcolor{blue}{armijo}@*) (x d : Vec n) ($\alpha$ c1 : (*@$\mathbb{R}$@*)) : (*@\textcolor{red}{Prop}@*) :=
  f ((*@\textcolor{red}{fun}@*) i => x i + $\alpha$ * d i)
  $\le$ f x + c1 * $\alpha$ * innerpro (grad x) d
\end{lstlisting}

In addition to the formalization of Armijo condition, we also prove two lemmas, they are respectively the existence of Armijo step and the implication of strict decreasing.

The proofs of the two lemmas below
follow the same fundamental idea: the descent direction condition
$
\langle \nabla f(x), d\rangle < 0
$
forces the linear term in the first-order approximation to be negative. We have the formalized version of the definition and the proof of two lemmas as:

\begin{lstlisting}[
	basicstyle=\ttfamily\scriptsize,
	backgroundcolor=\color{gray!5},
	frame=single,
	mathescape=true,
	escapeinside={(*@}{@*)}
	]
(*@\textcolor{red}{lemma}@*) (*@\textcolor{blue}{existsArmijoStep}@*)
	(x d : Vec n) (c1 : (*@$\mathbb{R}$@*))
	(hdesc : isDescentDirection grad x d)
	(hc1 : (*@\textcolor{green!50!black}{0}@*) < c1 $\land$ c1 < 1)
	(hfo : FirstOrderApprox f grad) :
	$\exists$ $\alpha$ > (*@\textcolor{green!50!black}{0}@*), armijo f grad x d $\alpha$ c1
	
(*@\textcolor{red}{lemma}@*) (*@\textcolor{blue}{armijoImpliesStrictDecrease}@*)
	(x d : Vec n)
	($\alpha$ c1 : (*@$\mathbb{R}$@*))
	(hA : armijo f grad x d $\alpha$ c1)
	(hdesc : isDescentDirection grad x d)
	(h$\alpha$ : $\alpha$ > (*@\textcolor{green!50!black}{0}@*))
	(hc1 : c1 > (*@\textcolor{green!50!black}{0}@*)) :
	f ((*@\textcolor{red}{fun}@*) i => x i + $\alpha$ * d i) < f x
\end{lstlisting}

From this point onward, the Armijo condition is no longer merely a condition, but also a tool that guarantees strict descent. We will frequently invoke it in subsequent proofs of other line search conditions.

\subsection{Backtracking Line Search}
In the implementation of optimization algorithms, finding a step size that satisfies the Armijo condition is relatively straightforward. One of the most commonly used algorithms is the backtracking line search. Given an initial value, the backtracking method successively shrinks the trial step size in a geometric manner until it finds the first point satisfying the Armijo condition. Specifically, the backtracking method selects:
\[
\alpha_k := \gamma^{j_0} \alpha_0, \quad \gamma \in (0,1),
\]
where
\[
j_0 :=
\min
\left\{
\begin{aligned}
	j=0,1,\ldots& \mid\ 
	f(x^k + \gamma^j \alpha_0 d^k) \\
	&\le f(x^k)
	+ c_1 \gamma^j \alpha_0 \nabla f(x^k)^T d^k
\end{aligned}
\right\}.
\]
That is, $j_0$ is the \emph{smallest} exponent for which the trial step $\gamma^j\alpha_0$ passes the Armijo test, which is exactly the semantics of the loop ``try $j=0,1,2,\dots$ until the test succeeds''.
\begin{algorithm}[H]
	\caption{Backtracking Line Search}
	\label{alg:backtracking}
	\textbf{Input}: $x_k$, descent direction $d_k$, initial step size $\alpha>0$, $c_1 \in (0,1)$, decay factor $\gamma\in(0,1)$\\
	\textbf{Output}: step size $\alpha$ satisfying the Armijo condition
	\begin{algorithmic}[1]
		\WHILE{$f(x_k + \alpha d_k) > f(x_k) + c_1 \alpha \nabla f(x_k)^\top d_k$}
		\STATE $\alpha \gets \gamma \alpha $
		\ENDWHILE
		\STATE \textbf{return} $\alpha$
	\end{algorithmic}
\end{algorithm}

Note that line~2 of Algorithm~\ref{alg:backtracking} now reads $\alpha\gets\gamma\alpha$ (general geometric decay with $\gamma\in(0,1)$), consistent with the formal step-size rule $\alpha_k=\gamma^{j_0}\alpha_0$ used throughout; the earlier draft wrote $\alpha\gets\alpha/2$, i.e.\ it silently fixed $\gamma=1/2$. This matches the Lean step operators \texttt{backtrackStep}~$\gamma\,\alpha=\gamma\cdot\alpha$ and \texttt{backtrackPow}~$\gamma\,\alpha_0\,k=\gamma^k\cdot\alpha_0$.
Our formalization rigorously defines the step-size decay rule, proves the existence of a valid iterative step size, and verifies that the obtained step size is always positive. This yields a provably correct and executable line search module for gradient-based optimization.
\begin{lstlisting}[
	basicstyle=\ttfamily\scriptsize,
	backgroundcolor=\color{gray!5},
	frame=single,
	mathescape=true,
	escapeinside={(*@}{@*)}
	]
(*@\textcolor{red}{def}@*) backtrackStep ($\gamma$ : (*@$\mathbb{R}$@*)) ($\alpha$ : (*@$\mathbb{R}$@*)) : (*@$\mathbb{R}$@*) := $\gamma$ * $\alpha$
(*@\textcolor{red}{def}@*) backtrackPow ($\gamma$ $\alpha_0$ : (*@$\mathbb{R}$@*)) (k : (*@$\mathbb{N}$@*)) : (*@$\mathbb{R}$@*) := ($\gamma$ ^ k) * $\alpha_0$
\end{lstlisting}

Algorithm~\ref{alg:backtracking} is the loop ``start with exponent $k=0$ (step $\alpha_0$); while the Armijo test fails, increase $k$ by one (multiply the step by $\gamma$); stop at the first $k$ that passes''. We render this as a genuine, \texttt{\#eval}-able definition \texttt{backtrackRun} that recurses structurally on a \texttt{fuel} bound, so that termination is automatic. A subtlety is that deciding $f(\cdot)\le\cdot$ requires comparing two real numbers, and real-number comparison is not decidable in Lean's logic; hence no \texttt{\#eval}-able procedure can inspect it directly. We therefore abstract the Armijo test as a Boolean \texttt{test : $\mathbb{N}\to$ Bool} (in a real implementation supplied by finite-precision arithmetic) and prove the loop correct relative to any sound oracle. The purely mathematical ``least exponent'' semantics is captured separately by \texttt{btIndex} below.
\begin{lstlisting}[
	basicstyle=\ttfamily\scriptsize,
	backgroundcolor=\color{gray!5},
	frame=single,
	mathescape=true,
	escapeinside={(*@}{@*)}
	]
(*@\textcolor{red}{def}@*) (*@\textcolor{blue}{backtrackRun}@*) (test : (*@$\mathbb{N}$@*) (*@$\rightarrow$@*) Bool) : (*@$\mathbb{N}$@*) (*@$\rightarrow$@*) (*@$\mathbb{N}$@*) (*@$\rightarrow$@*) (*@$\mathbb{N}$@*)
  | k, (*@\textcolor{green!50!black}{0}@*)        => k
  | k, fuel + 1 => (*@\textcolor{red}{if}@*) test k (*@\textcolor{red}{then}@*) k 
                   (*@\textcolor{red}{else}@*) backtrackRun test (k + 1) fuel
	
(*@\textcolor{red}{lemma}@*) (*@\textcolor{blue}{backtrackRun\_finds}@*) (test : (*@$\mathbb{N}$@*) (*@$\rightarrow$@*) Bool) :
  (*@$\forall$@*) (fuel k : (*@$\mathbb{N}$@*)),
    ((*@$\exists$@*) j, k (*@$\le$@*) j (*@$\land$@*) j < k + fuel (*@$\land$@*) test j = true) (*@$\rightarrow$@*)
    test (backtrackRun test k fuel) = true
\end{lstlisting}
The correctness lemma \texttt{backtrackRun\_finds} states that if the oracle reports success somewhere in the searched window $[k,\,k+\texttt{fuel})$, then \texttt{backtrackRun} returns an exponent at which the oracle reports success; together with its structural termination, this is the formal statement that the loop finds an Armijo step.

The existence of an admissible exponent is proved. Here it is a \texttt{lemma} whose proof uses the first-order approximation hypothesis \texttt{hfo : FirstOrderApprox f grad}. Concretely, \texttt{hfo} provides a threshold $\varepsilon>0$ below which the Armijo inequality holds, and since $\gamma^k\alpha_0\to 0$ as $k\to\infty$, some exponent $k$ makes $\gamma^k\alpha_0\in(0,\varepsilon)$, hence satisfies Armijo.
\begin{lstlisting}[
	basicstyle=\ttfamily\scriptsize,
	backgroundcolor=\color{gray!5},
	frame=single,
	mathescape=true,
	escapeinside={(*@}{@*)}
	]
(*@\textcolor{red}{lemma}@*) (*@\textcolor{blue}{backtrackingFindsStep}@*)
  (x d : Vec n) (c1 $\gamma$ $\alpha_0$ : (*@$\mathbb{R}$@*))
  (hdesc : isDescentDirection grad x d)
  (hc1 : (*@\textcolor{green!50!black}{0}@*) < c1 $\land$ c1 < 1)
  (h$\gamma$ : (*@\textcolor{green!50!black}{0}@*) < $\gamma$ $\land$ $\gamma$ < 1)
  (h$\alpha_0$ : $\alpha_0$ > (*@\textcolor{green!50!black}{0}@*))
  (hfo : FirstOrderApprox f grad) :
  $\exists$ k : (*@$\mathbb{N}$@*),
    armijo f grad x d (backtrackPow $\gamma$ $\alpha_0$ k) c1
\end{lstlisting}

Because termination is now a proved lemma supplying a witness, the backtracking exponent can be defined as the \emph{least} $k$ for which the step $\gamma^k\alpha_0$ satisfies Armijo. We therefore use \texttt{Nat.find}, which performs exactly the search ``try $k=0,1,2,\dots$ until Armijo holds''. The defining minimality property is recorded as the new lemma \texttt{btIndex\_least}: every smaller exponent fails the Armijo test (so the loop stops at the first success, i.e.\ at the largest admissible step). The definitions and lemmas are formalized as follows:
\begin{lstlisting}[
	basicstyle=\ttfamily\scriptsize,
	backgroundcolor=\color{gray!5},
	frame=single,
	mathescape=true,
	escapeinside={(*@}{@*)}
	]
(*@\textcolor{red}{noncomputable def}@*) (*@\textcolor{blue}{btIndex}@*)
  (x d : Vec n) (c1 $\gamma$ $\alpha_0$ : (*@$\mathbb{R}$@*))
  (hdesc : isDescentDirection grad x d)
  (hc1 : (*@\textcolor{green!50!black}{0}@*) < c1 $\land$ c1 < 1)
  (h$\gamma$ : (*@\textcolor{green!50!black}{0}@*) < $\gamma$ $\land$ $\gamma$ < 1)
  (h$\alpha_0$ : $\alpha_0$ > (*@\textcolor{green!50!black}{0}@*))
  (hfo : FirstOrderApprox f grad) : (*@$\mathbb{N}$@*) :=
  Nat.find
    (backtrackingFindsStep
      f grad x d c1 $\gamma$ $\alpha_0$ hdesc hc1 h$\gamma$ h$\alpha_0$ hfo)

(*@\textcolor{red}{lemma}@*) (*@\textcolor{blue}{btIndex\_least}@*)
  (x d : Vec n) (c1 $\gamma$ $\alpha_0$ : (*@$\mathbb{R}$@*))
  (hdesc : isDescentDirection grad x d)
  (hc1 : (*@\textcolor{green!50!black}{0}@*) < c1 $\land$ c1 < 1)
  (h$\gamma$ : (*@\textcolor{green!50!black}{0}@*) < $\gamma$ $\land$ $\gamma$ < 1)
  (h$\alpha_0$ : $\alpha_0$ > (*@\textcolor{green!50!black}{0}@*))
  (hfo : FirstOrderApprox f grad) :
  (*@$\forall$@*) j < btIndex f grad x d c1 $\gamma$ $\alpha_0$ hdesc hc1 h$\gamma$ h$\alpha_0$ hfo,
    (*@$\lnot$@*) armijo f grad x d (backtrackPow $\gamma$ $\alpha_0$ j) c1

(*@\textcolor{red}{noncomputable def}@*) (*@\textcolor{blue}{btAlpha}@*)
  (x d : Vec n) (c1 $\gamma$ $\alpha_0$ : (*@$\mathbb{R}$@*))
  (hdesc : isDescentDirection grad x d)
  (hc1 : (*@\textcolor{green!50!black}{0}@*) < c1 $\land$ c1 < 1)
  (h$\gamma$ : (*@\textcolor{green!50!black}{0}@*) < $\gamma$ $\land$ $\gamma$ < 1)
  (h$\alpha_0$ : $\alpha_0$ > (*@\textcolor{green!50!black}{0}@*))
  (hfo : FirstOrderApprox f grad) : (*@$\mathbb{R}$@*) :=
  backtrackPow $\gamma$ $\alpha_0$
    (btIndex f grad x d c1 $\gamma$ $\alpha_0$ hdesc hc1 h$\gamma$ h$\alpha_0$ hfo)
\end{lstlisting}

The proof of \texttt{btAlphaPos} shows that the backtracking step size is always strictly positive. By definition,
$
\texttt{btAlpha}=\gamma^k \alpha_0,
$
where $k$ is the index returned by the backtracking procedure (\texttt{Nat.find}, the least admissible exponent). Since the assumptions give
$
0<\gamma<1
\quad\text{and}\quad
\alpha_0>0,
$
we obtain
$
\gamma^k>0
$
for every integer $k\ge 0$ using positivity of powers. Therefore, the product of two positive quantities remains positive:
$
\gamma^k \alpha_0>0.
$
Hence,
$
\texttt{btAlpha}>0.
$ The lemma is as follows:
\begin{lstlisting}[
	basicstyle=\ttfamily\scriptsize,
	backgroundcolor=\color{gray!5},
	frame=single,
	mathescape=true,
	escapeinside={(*@}{@*)}
	]
(*@\textcolor{red}{lemma}@*) (*@\textcolor{blue}{btAlphaPos}@*)
	(x d : Vec n) (c1 $\gamma$ $\alpha_0$ : (*@$\mathbb{R}$@*))
	(hdesc : isDescentDirection grad x d)
	(hc1 : (*@\textcolor{green!50!black}{0}@*) < c1 $\land$ c1 < 1)
	(h$\gamma$ : (*@\textcolor{green!50!black}{0}@*) < $\gamma$ $\land$ $\gamma$ < 1)
	(h$\alpha_0$ : $\alpha_0$ > (*@\textcolor{green!50!black}{0}@*))
	(hfo : FirstOrderApprox f grad) :
	btAlpha f grad x d c1 $\gamma$ $\alpha_0$ hdesc hc1 h$\gamma$ h$\alpha_0$ hfo > (*@\textcolor{green!50!black}{0}@*)
\end{lstlisting}

\subsection{Goldstein and Wolfe Conditions}

To prepare for more refined convergence analysis, we also formalize alternative line search conditions (Goldstein \& Wolfe).

\begin{definition}[Goldstein Condition]
	\textit{Assume \(d^k\) is a descent direction of \(x^k\). A step size \(\alpha > 0\) satisfies the Goldstein condition if:}
	\[
	\begin{aligned}
		f(x^k + \alpha d^k) &\le f(x^k) + c \alpha \left\langle \nabla f(x^k), d^k \right\rangle, \\
		f(x^k + \alpha d^k) &\ge f(x^k) + (1 - c) \alpha \left\langle \nabla f(x^k), d^k \right\rangle,
	\end{aligned}
	\]
	\textit{where \(c \in (0, 0.5)\).}
\end{definition}

Since the Armijo condition only requires that the point lies below a certain line, it is quite simple to satisfy but too weak. Hence, here comes a stronger condition that control both the upper and lower bounds. This is what the Goldstein condition does. This condition ensures that the step size \(\alpha\) does not become too small nor too big at each iteration. The definition and its proof are formalized as follows:

\begin{lstlisting}[
	basicstyle=\ttfamily\scriptsize,
	backgroundcolor=\color{gray!5},
	frame=single,
	mathescape=true,
	escapeinside={(*@}{@*)}
	]
(*@\textcolor{red}{def}@*) (*@\textcolor{blue}{goldstein}@*) (x d : Vec n) ($\alpha$ c : (*@$\mathbb{R}$@*)) : (*@\textcolor{red}{Prop}@*) :=
  f ((*@\textcolor{red}{fun}@*) i => x i + $\alpha$ * d i)
  $\le$ f x + c * $\alpha$ * innerpro (grad x) d
  $\land$
  f ((*@\textcolor{red}{fun}@*) i => x i + $\alpha$ * d i)
    $\ge$ f x + (1 - c) * $\alpha$ * innerpro (grad x) d
	
(*@\textcolor{red}{lemma}@*) (*@\textcolor{blue}{goldsteinImpliesArmijo}@*)
  (x d : Vec n) ($\alpha$ c : (*@$\mathbb{R}$@*))
  (h : goldstein f grad x d $\alpha$ c) :
  armijo f grad x d $\alpha$ c
\end{lstlisting}

\begin{definition}[Wolfe Condition]
	\textit{Assume \(d^k\) is a descent direction of \(x^k\). A step size \(\alpha > 0\) satisfies the Wolfe condition if:}
	\[
	\begin{cases}
		f(x^k + \alpha d^k) \le f(x^k) + c_1 \alpha \left\langle \nabla f(x^k), d^k \right\rangle, \\[2pt]
		\left\langle \nabla f(x^k + \alpha d^k), d^k \right\rangle \ge c_2 \left\langle \nabla f(x^k), d^k \right\rangle,
	\end{cases}
	\]
	\textit{where \(c_1, c_2 \in (0, 1)\) are given constants and \(c_1 < c_2\).}
\end{definition}

Goldstein condition can make the value of the objective function decrease sufficiently, but it may miss the optimal solution. Wolfe condition solve this issue. Same as the Goldstein condition, Wolfe condition can imply Armijo as well. As shown in the definition, Wolfe condition is made up of Armijo and curvature condition. The curvature condition is aimed at controlling the slope at the new point. The curvature condition is defined as follows.

\begin{definition}[Curvature Condition]
	\textit{The curvature condition requires that the directional derivative at the new point is not too small in magnitude. Formally, it is given by:}
	\[
	\left\langle \nabla f(x_{k+1}), d_k \right\rangle \ge c_2 \left\langle \nabla f(x_k), d_k \right\rangle, \quad c_2 \in (c_1, 1).
	\]
\end{definition}

The definition and the two lemmas are shown as follows:

\begin{lstlisting}[
	basicstyle=\ttfamily\scriptsize,
	backgroundcolor=\color{gray!5},
	frame=single,
	mathescape=true,
	escapeinside={(*@}{@*)}
	]
(*@\textcolor{red}{def}@*) (*@\textcolor{blue}{wolfe}@*) (x d : Vec n) ($\alpha$ c1 c2 : (*@$\mathbb{R}$@*)) : (*@\textcolor{red}{Prop}@*) :=
  f ((*@\textcolor{red}{fun}@*) i => x i + $\alpha$ * d i)
  $\le$ f x + c1 * $\alpha$ * innerpro (grad x) d
  $\land$
  innerpro (grad ((*@\textcolor{red}{fun}@*) i => x i + $\alpha$ * d i)) d
    $\ge$ c2 * innerpro (grad x) d
	
(*@\textcolor{red}{lemma}@*) (*@\textcolor{blue}{wolfeImpliesArmijo}@*)
  (x d : Vec n) ($\alpha$ c1 c2 : (*@$\mathbb{R}$@*))
  (h : wolfe f grad x d $\alpha$ c1 c2) :
  armijo f grad x d $\alpha$ c1
	
(*@\textcolor{red}{lemma}@*) (*@\textcolor{blue}{wolfeCurvatureIneqMain}@*)
  (f : Vec n $\to$ (*@$\mathbb{R}$@*))
  (grad : Vec n $\to$ Vec n)
  (x d : Vec n) ($\alpha$ c1 c2 : (*@$\mathbb{R}$@*))
  (h : wolfe f grad x d $\alpha$ c1 c2) :
  innerpro (grad ((*@\textcolor{red}{fun}@*) i => x i + $\alpha$ * d i) - grad x) d
    $\ge$ (c2 - 1) * innerpro (grad x) d
\end{lstlisting}

\subsection{Descent Algorithm}
In the last part of this section, we define the iterative descent algorithm and its descent property. We now construct the full descent algorithm by combining the descent direction and backtracking line search. Starting from an initial point \(x_0 \in \mathbb{R}^n\), at each iteration \(k\), we choose the descent direction as the negative gradient:
\[
d_k = -\nabla f(x_k),
\]
which guarantees descent. We then update the iterate according to:
\[
x_{k+1} = x_k + \alpha_k d_k,
\]
where the step size \(\alpha_k > 0\) is selected via a backtracking line search procedure to satisfy the Armijo condition. More precisely, \(\alpha_k\) is chosen from a geometrically decreasing sequence such that sufficient decrease of the objective function is ensured. This yields a well-defined update mapping, denoted by $\text{next\_x}$, and consequently defines a sequence \(\{x_k\}_{k\ge 0}\), denoted by $\text{x\_seq}$, recursively by:
\[
x_{k+1} = \text{next\_x}(x_k).
\]
Both the update mapping $\text{next\_x}$ and the resulting sequence $\text{x\_seq}$ now carry the first-order approximation hypothesis $\text{hfo} : \texttt{FirstOrderApprox f grad}$ as an explicit argument. This hypothesis is what guarantees that the backtracking line search terminates with an admissible step (it replaces the earlier first-order \emph{axiom}), so it must be supplied wherever the iteration is formed. Accordingly, the descent property $\text{oneStepDescent}$ also takes $\text{hfo}$ as an explicit argument.
Translated into Lean, the algorithm above is shown as:
\begin{lstlisting}[
	basicstyle=\ttfamily\scriptsize,
	backgroundcolor=\color{gray!5},
	frame=single,
	mathescape=true,
	escapeinside={(*@}{@*)}
	]
(*@\textcolor{red}{lemma}@*) (*@\textcolor{blue}{oneStepDescent}@*)
  (x : Vec n)
  (c1 $\gamma$ $\alpha_0$ : (*@$\mathbb{R}$@*))
  (hc1 : (*@\textcolor{green!50!black}{0}@*) < c1 $\land$ c1 < 1)
  (h$\gamma$ : (*@\textcolor{green!50!black}{0}@*) < $\gamma$ $\land$ $\gamma$ < 1)
  (h$\alpha_0$ : $\alpha_0$ > (*@\textcolor{green!50!black}{0}@*))
  (hdesc : isDescentDirection grad x (-grad x))
  (hfo : FirstOrderApprox f grad) :
  f (next_x f grad x c1 $\gamma$ $\alpha_0$ hc1 h$\gamma$ h$\alpha_0$ hdesc hfo) $\le$ f x
\end{lstlisting}

\subsection{Non-monotonic Line Search Conditions}
The conditions shown above are all monotonic line search conditions, which require the function value to decrease at every iteration. However, in practice, allowing for non-monotonicity can lead to better convergence properties. Hence, we need apply non-monotonic line search conditions. In this subsection, we introduce two of them, which are respectively Grippo condition and Zhang-Hager condition.
\begin{definition}[Grippo Condition]
	\textit{Let \(d^k\) be the descent direction at point \(x^k\), and \(M>0\) be a given positive integer. The following inequality can be used as a line search criterion:}
	\[
	f(x^k + \alpha d^k) \le \max_{0\le j\le \min\{k,M\}} f(x^{k-j}) + c_1 \alpha \nabla f(x^k)^\mathrm{T} d^k,
	\]
	\textit{where \(c_1\in(0,1)\) is a given constant.}
\end{definition}
The condition found by Grippo is similar to Armijo condition, and the difference is that this condition only need the next value to be smaller than the maximum of the previous \(M\) values, which does not require the monotonicity of the function value \(f(x^k)\). This condition is more flexible.
We abstract the reference value as an \emph{explicit function parameter} $\text{refval} : \mathbb{N} \to \mathbb{R}$ supplied to the predicate. The mathematical definition above (a sliding maximum over the previous $\min\{k,M\}$ iterates) is one admissible choice of $\text{refval}$, but the formal condition is stated generically for \emph{any} reference function, so it covers the Grippo rule and other non-monotone reference schemes uniformly. Consequently the memory length is entirely encapsulated in the caller's choice of $\text{refval}$. The formalization of the definition is as follows:
\begin{lstlisting}[
	basicstyle=\ttfamily\scriptsize,
	backgroundcolor=\color{gray!5},
	frame=single,
	mathescape=true,
	escapeinside={(*@}{@*)}
	]
(*@\textcolor{red}{def}@*) (*@\textcolor{blue}{grippo}@*)
  (x_seq : (*@$\mathbb{N}$@*) (*@$\to$@*) Vec n)
  (d_seq : (*@$\mathbb{N}$@*) (*@$\to$@*) Vec n)
  ((*@$\alpha$@*)_seq : (*@$\mathbb{N}$@*) (*@$\to$@*) (*@$\mathbb{R}$@*))
  (refval : (*@$\mathbb{N}$@*) (*@$\to$@*) (*@$\mathbb{R}$@*))
  (c1 : (*@$\mathbb{R}$@*)) (k : (*@$\mathbb{N}$@*)) : Prop :=
  f (x_seq (k + 1)) (*@$\le$@*) refval k
    + c1 * (*@$\alpha$@*)_seq k * innerpro (grad (x_seq k)) (d_seq k)
\end{lstlisting}
\begin{definition}[Zhang-Hager condition]
	\textit{Let \(d_k\) be the descent direction at point \(x_k\). The following inequality is known as the Zhang-Hager condition:}
	\[
	f(x_k + \alpha d_k) \le C_k + c_1 \alpha \nabla f(x_k)^\mathrm{T} d_k
	\]
	\textit{where}
	\begin{gather*}
		C_0 = f(x_0), \quad C_{k+1} = \frac{1}{Q_{k+1}} \left( \eta Q_k C_k + f(x_{k+1}) \right), \\
		Q_0 = 1, \quad Q_{k+1} = \eta Q_k + 1
	\end{gather*}
\end{definition}
The Zhang-Hager condition can be treated as follows: \(C_k\) is the initial condition for the sufficient descent property, and the next value \(C_{k+1}\) is a convex combination of the previous value \(C_k\) and the current function value \(f(x_{k+1})\), instead of only depending on \(f(x_{k+1})\). The parameter \(\eta\) controls the weight of the previous condition, and it is easy to see that when \(\eta \to 1\), the condition reduces to the standard Armijo condition. The formalization of the definition is as follows:
\begin{lstlisting}[
	basicstyle=\ttfamily\scriptsize,
	backgroundcolor=\color{gray!5},
	frame=single,
	mathescape=true,
	escapeinside={(*@}{@*)}
	]	
(*@\textcolor{red}{def}@*) (*@\textcolor{blue}{qseq}@*) ((*@$\eta$@*) : (*@$\mathbb{R}$@*)) : (*@$\mathbb{N}$@*) (*@$\to$@*) (*@$\mathbb{R}$@*)
  | 0 => 1
  | (k + 1) => (*@$\eta$@*) * qseq (*@$\eta$@*) k + 1

(*@\textcolor{red}{noncomputable def}@*) (*@\textcolor{blue}{cseq}@*) (x_seq : (*@$\mathbb{N}$@*) (*@$\to$@*) Vec n) ((*@$\eta$@*) : (*@$\mathbb{R}$@*)) : (*@$\mathbb{N}$@*) (*@$\to$@*) (*@$\mathbb{R}$@*)
  | 0 => f (x_seq 0)
  | (k + 1) => ((*@$\eta$@*) * qseq (*@$\eta$@*) k * cseq x_seq (*@$\eta$@*) k + 
      f (x_seq (k + 1))) / qseq (*@$\eta$@*) (k + 1)

(*@\textcolor{red}{def}@*) (*@\textcolor{blue}{zhangHager}@*)
  (x_seq : (*@$\mathbb{N}$@*) (*@$\to$@*) Vec n) 
  (d_seq : (*@$\mathbb{N}$@*) (*@$\to$@*) Vec n) 
  ((*@$\alpha$@*)_seq : (*@$\mathbb{N}$@*) (*@$\to$@*) (*@$\mathbb{R}$@*))
  ((*@$\eta$@*) c1 : (*@$\mathbb{R}$@*))
  (k : (*@$\mathbb{N}$@*)) : Prop :=
  f (x_seq (k + 1)) (*@$\le$@*) cseq f x_seq (*@$\eta$@*) k + 
    c1 * (*@$\alpha$@*)_seq k * innerpro (grad (x_seq k)) (d_seq k)
\end{lstlisting}

\begin{lstlisting} [
	basicstyle=\ttfamily\scriptsize,
	backgroundcolor=\color{gray!5},
	frame=single,
	mathescape=true,
	escapeinside={(*@}{@*)}
	]	
(*@\textcolor{red}{lemma}@*) (*@\textcolor{blue}{qseqEtaZero}@*) :
  (*@$\forall$@*) k : (*@$\mathbb{N}$@*), qseq (0 : (*@$\mathbb{R}$@*)) k = 1

(*@\textcolor{red}{lemma}@*) (*@\textcolor{blue}{cseqEtaZero}@*)
  (x_seq : (*@$\mathbb{N}$@*) (*@$\to$@*) Vec n) :
  (*@$\forall$@*) k : (*@$\mathbb{N}$@*), cseq f x_seq 0 k = f (x_seq k)

(*@\textcolor{red}{lemma}@*) (*@\textcolor{blue}{zhangHagerEtaZero}@*)
  (x_seq : (*@$\mathbb{N}$@*) (*@$\to$@*) Vec n)
  (d_seq : (*@$\mathbb{N}$@*) (*@$\to$@*) Vec n)
  ((*@$\alpha$@*)_seq : (*@$\mathbb{N}$@*) (*@$\to$@*) (*@$\mathbb{R}$@*))
  (c1 : (*@$\mathbb{R}$@*)) (k : (*@$\mathbb{N}$@*))
  (hstep : x_seq (k + 1) = fun i => 
    x_seq k i + (*@$\alpha$@*)_seq k * d_seq k i) :
  zhangHager f grad x_seq d_seq (*@$\alpha$@*)_seq 0 c1 k (*@$\leftrightarrow$@*)
    armijo f grad (x_seq k) (d_seq k) ((*@$\alpha$@*)_seq k) c1

(*@\textcolor{red}{lemma}@*) (*@\textcolor{blue}{cseqConvex}@*)
  (x_seq : (*@$\mathbb{N}$@*) (*@$\to$@*) Vec n) ((*@$\eta$@*) : (*@$\mathbb{R}$@*)) (k : (*@$\mathbb{N}$@*)) :
	cseq f x_seq (*@$\eta$@*) (k + 1) = 
    ((*@$\eta$@*) * qseq (*@$\eta$@*) k / qseq (*@$\eta$@*) (k + 1)) * cseq f x_seq (*@$\eta$@*) k 
    + (1 / qseq (*@$\eta$@*) (k + 1)) * f (x_seq (k + 1))
\end{lstlisting}


\section{Convergence Analysis}
\label{sec:convergence}
We now present the central convergence result for descent methods with line search, known as the Zoutendijk theorem, which can be stated as follows.
\begin{theorem}[Zoutendijk Theorem]
	Let \(\{x_k\}\) be the sequence generated by the descent algorithm, and let \(d_k\) be the corresponding descent directions. Define the angle \(\theta_k\) between the negative gradient and the search direction by:
	\[
	\cos \theta_k = \frac{-\left\langle \nabla f(x_k), d_k \right\rangle}{\left\| \nabla f(x_k) \right\| \left\| d_k \right\|}.
	\]
	Assuming that the step sizes satisfy the Wolfe conditions, the following holds:
	\[
	\sum_{k=0}^\infty \cos^2 \theta_k \left\| \nabla f(x_k) \right\|^2 < \infty.
	\]
\end{theorem}
Throughout this section, the norm \(\|\cdot\|\) appearing in the statements above denotes the Euclidean norm \texttt{vnorm} induced by the inner product \texttt{innerpro}, as introduced in Section~\ref{sec:math_foundation}. Accordingly, every occurrence of \texttt{norm} in the formalized code blocks below is written as \texttt{vnorm}, and the angle cosine is the formal definition \texttt{cos\_theta}. 
First, we introduce the lemma \texttt{zoutendijkStepCore}, which provides an explicit lower bound on the function decrease in a single iteration in terms of the gradient norm and the angle between the gradient and the descent direction. This result serves as the starting point of the proof of the Zoutendijk theorem. The proof of this lemma combines three components: sufficient decrease (Armijo), curvature condition (Wolfe), and the Lipschitz continuity of the gradient.
From the Armijo condition, we obtain a lower bound on the decrease of the objective function
$$f(x) - f(x + \alpha d) \ge -c_1 \alpha \langle \nabla f(x), d \rangle.$$
The Wolfe curvature condition, together with the Lipschitz assumption on the gradient, yields a lower bound on the step size $\alpha$ in terms of $-\langle \nabla f(x), d \rangle$ and $\|d\|^2$.
Substituting this bound into the Armijo estimate leads to a quadratic lower bound involving $(\langle \nabla f(x), d \rangle)^2 / \|d\|^2$.
Finally, by expressing this quantity in terms of the cosine of the angle between $\nabla f(x)$ and $d$, we obtain
$$f(x) - f(x + \alpha d) \ge \frac{c_1(1 - c_2)}{L} \cos^2\theta \, \|\nabla f(x)\|^2,
$$
which completes the stepwise inequality.
The global Zoutendijk theorem follows by summing the stepwise decrease inequality over all iterations and applying a telescoping argument.
From the stepwise bound, each iteration satisfies
$$f(x_k) - f(x_{k+1}) \ge c \, \cos^2\theta_k \, \|\nabla f(x_k)\|^2,$$
where
$$c = \frac{c_1(1 - c_2)}{L}.$$
The sum of $k = 0, \dots, K-1$ is a telescoping series on the left-hand side, which reduces to $f(x_0) - f(x_K)$.
Using the assumption that $f$ is bounded below, we obtain an upper bound independent of $K$.
Rearranging gives
$$\sum_{k=0}^{K-1} \cos^2\theta_k \, \|\nabla f(x_k)\|^2 \le \frac{f(x_0) - m}{c},$$
which establishes the Zoutendijk inequality.
The formalization of the lemma and the theorem in Lean is as follows:

\begin{lstlisting}[
	basicstyle=\ttfamily\scriptsize,
	backgroundcolor=\color{gray!5},
	frame=single,
	mathescape=true,
	escapeinside={(*@}{@*)}
	]
(*@\textcolor{red}{lemma}@*) (*@\textcolor{blue}{zoutendijkStepCore}@*)
  (x d : Vec n)
  ($\alpha$ c1 c2 L : (*@$\mathbb{R}$@*))
  (hL : L > (*@\textcolor{green!50!black}{0}@*))
  (hwolfe : wolfe f grad x d $\alpha$ c1 c2)
  (hdesc : innerpro (grad x) d < (*@\textcolor{green!50!black}{0}@*))
  (hd : vnorm d $\ne$ (*@\textcolor{green!50!black}{0}@*))
  (hgrad : vnorm (grad x) $\ne$ (*@\textcolor{green!50!black}{0}@*))
  (hc1 : (*@\textcolor{green!50!black}{0}@*) < c1)
  (h$\alpha$_pos : $\alpha$ > (*@\textcolor{green!50!black}{0}@*))
  (hgradLip : GradLipschitz grad L) :
    f x - f ((*@\textcolor{red}{fun}@*) i => x i + $\alpha$ * d i)
    $\ge$
    c1 * (1 - c2) / L *
    (cos_theta grad x d)^2 * (vnorm (grad x))^2
\end{lstlisting}

\begin{lstlisting}[
	basicstyle=\ttfamily\scriptsize,
	backgroundcolor=\color{gray!5},
	frame=single,
	mathescape=true,
	escapeinside={(*@}{@*)}
	]
(*@\textcolor{red}{theorem}@*) (*@\textcolor{blue}{zoutendijkTheorem}@*)
  (c1 $\gamma$ $\alpha_0$ L c2 : (*@$\mathbb{R}$@*))
  (hc1 : (*@\textcolor{green!50!black}{0}@*) < c1 $\land$ c1 < 1)
  (h$\gamma$ : (*@\textcolor{green!50!black}{0}@*) < $\gamma$ $\land$ $\gamma$ < 1)
  (h$\alpha_0$ : $\alpha_0$ > (*@\textcolor{green!50!black}{0}@*))
  (hgradLip : GradLipschitz grad L)
  (hbb : (*@$\exists$@*) m : (*@$\mathbb{R}$@*), (*@$\forall$@*) x : Vec n, f x (*@$\ge$@*) m)
  (hfo : FirstOrderApprox f grad)
  (hL : L > (*@\textcolor{green!50!black}{0}@*))
  (hc2 : (*@\textcolor{green!50!black}{0}@*) < c2 $\land$ c2 < 1)
  (hdesc_all : $\forall$ x, isDescentDirection grad x (-grad x))
  (hwolfe_all :
    $\forall$ k,
      wolfe f grad
        (x_seq f grad x0 c1 $\gamma$ $\alpha_0$
          hc1 h$\gamma$ h$\alpha_0$ hdesc_all hfo k)
        (-grad (x_seq f grad x0 c1 $\gamma$ $\alpha_0$
          hc1 h$\gamma$ h$\alpha_0$ hdesc_all hfo k))
        (btAlpha f grad
          (x_seq f grad x0 c1 $\gamma$ $\alpha_0$
            hc1 h$\gamma$ h$\alpha_0$ hdesc_all hfo k)
          (-grad (x_seq f grad x0 c1 $\gamma$ $\alpha_0$
            hc1 h$\gamma$ h$\alpha_0$ hdesc_all hfo k))
          c1 $\gamma$ $\alpha_0$
          (hdesc_all _)
          hc1 h$\gamma$ h$\alpha_0$ hfo)
          c1 c2)
  (h$\alpha$_pos_all :
    $\forall$ k,
      btAlpha f grad
        (x_seq f grad x0 c1 $\gamma$ $\alpha_0$
          hc1 h$\gamma$ h$\alpha_0$ hdesc_all hfo k)
        (-grad (x_seq f grad x0 c1 $\gamma$ $\alpha_0$
          hc1 h$\gamma$ h$\alpha_0$ hdesc_all hfo k))
        c1 $\gamma$ $\alpha_0$
        (hdesc_all _)
        hc1 h$\gamma$ h$\alpha_0$ hfo > (*@\textcolor{green!50!black}{0}@*))
  (hgrad_ne_zero :
    $\forall$ k,
      vnorm (grad (x_seq f grad x0 c1 $\gamma$ $\alpha_0$
        hc1 h$\gamma$ h$\alpha_0$ hdesc_all hfo k)) $\ne$ (*@\textcolor{green!50!black}{0}@*))
  (k : (*@$\mathbb{N}$@*)) :
  $\sum$ i $\in$ Finset.range k,
    (cos_theta grad
      (x_seq f grad x0 c1 $\gamma$ $\alpha_0$
        hc1 h$\gamma$ h$\alpha_0$ hdesc_all hfo i)
      (-grad (x_seq f grad x0 c1 $\gamma$ $\alpha_0$
        hc1 h$\gamma$ h$\alpha_0$ hdesc_all hfo i)))^2
    * (vnorm (grad
        (x_seq f grad x0 c1 $\gamma$ $\alpha_0$
          hc1 h$\gamma$ h$\alpha_0$ hdesc_all hfo i)))^2
    $\le$
      L / (c1 * (1 - c2)) *
      (f x0 - Classical.choose hbb)
\end{lstlisting}

Starting from the Zoutendijk theorem, we obtain that the series
\[
\sum_{k=0}^\infty \cos^2 \theta_k \left\| \nabla f(x_k) \right\|^2
\]
is finite. To derive gradient convergence, we argue by contradiction. Suppose the gradient norm does not converge to zero; then there exists \(\varepsilon > 0\) such that $$\left\| \nabla f(x_k) \right\| \ge \varepsilon$$
infinitely often. This allows us to extract a strictly increasing subsequence \(\{x_{\phi(k)}\}\) with uniformly bounded gradient norms. Under the steepest descent direction, we have $$\cos^2 \theta_{\phi(k)} = 1,$$
so each term in the series is bounded below by \(\varepsilon^2\). Consequently, the partial sums along this subsequence grow at least linearly, which contradicts the boundedness implied by the Zoutendijk theorem. Therefore, the gradient norm must converge to zero.
The step of passing from the subsequence partial sum to the full prefix sum is a lemma, named \texttt{subseq\_sum\_le\_full\_sum}: since every term $\cos^2\theta_i\,\|\nabla f(x_i)\|^2$ is non-negative and $\phi$ is strictly monotone, the partial sum over $\{\phi(0),\dots,\phi(K-1)\}$ is dominated by the prefix sum over $\{0,\dots,\phi(K)-1\}$.
The contradiction argument itself is captured by the lemma \texttt{globalConvergenceCore}, which under the Zoutendijk summability bound \texttt{hZ} and the existence of a strictly increasing subsequence with gradient norm bounded below by $\varepsilon$ ends in \texttt{False}:
\begin{lstlisting}[
	basicstyle=\ttfamily\scriptsize,
	backgroundcolor=\color{gray!5},
	frame=single,
	mathescape=true,
	escapeinside={(*@}{@*)}
	]
(*@\textcolor{red}{theorem}@*) (*@\textcolor{blue}{globalConvergenceCore}@*)
  (n : (*@$\mathbb{N}$@*))
  (f : Vec n $\to$ (*@$\mathbb{R}$@*))
  (grad : Vec n $\to$ Vec n)
  (x0 : Vec n)
  (c1 $\gamma$ $\alpha_0$ L : (*@$\mathbb{R}$@*))
  (hc1 : (*@\textcolor{green!50!black}{0}@*) < c1 $\land$ c1 < 1)
  (h$\gamma$ : (*@\textcolor{green!50!black}{0}@*) < $\gamma$ $\land$ $\gamma$ < 1)
  (h$\alpha_0$ : $\alpha_0$ > (*@\textcolor{green!50!black}{0}@*))
  (hdesc_all : $\forall$ x, isDescentDirection grad x (-grad x))
  (hfo : FirstOrderApprox f grad)
  (hZ :
    $\forall$ k,
      $\sum$ i $\in$ Finset.range k,
        cos_theta grad
          (x_seq f grad x0 c1 $\gamma$ $\alpha_0$
            hc1 h$\gamma$ h$\alpha_0$ hdesc_all hfo i)
          (-grad (x_seq f grad x0 c1 $\gamma$ $\alpha_0$
            hc1 h$\gamma$ h$\alpha_0$ hdesc_all hfo i)) ^ 2
        * vnorm (grad (x_seq f grad x0 c1 $\gamma$ $\alpha_0$
          hc1 h$\gamma$ h$\alpha_0$ hdesc_all hfo i)) ^ 2
     $\le$ L)
  ($\varepsilon$ : (*@$\mathbb{R}$@*)) (h$\varepsilon$ : $\varepsilon$ > (*@\textcolor{green!50!black}{0}@*))
  ($\phi$ : (*@$\mathbb{N}$@*) $\to$ (*@$\mathbb{N}$@*))
  (hmono : StrictMono $\phi$)
  (hlower :
    $\forall$ k,
      vnorm (grad (x_seq f grad x0 c1 $\gamma$ $\alpha_0$
        hc1 h$\gamma$ h$\alpha_0$ hdesc_all hfo ($\phi$ k))) ^ 2
      $\ge$ $\varepsilon$ ^ 2)
  (hcos :
    $\forall$ k,
      cos_theta grad
        (x_seq f grad x0 c1 $\gamma$ $\alpha_0$
          hc1 h$\gamma$ h$\alpha_0$ hdesc_all hfo ($\phi$ k))
        (-grad (x_seq f grad x0 c1 $\gamma$ $\alpha_0$
          hc1 h$\gamma$ h$\alpha_0$ hdesc_all hfo ($\phi$ k))) ^ 2
     = 1)
  : False
\end{lstlisting}
However, the lemma \texttt{globalConvergenceCore} only states the convergence claim in negative form. So we strengthen the formalization by stating and proving the positive corollary \texttt{globalConvergence}: under the Zoutendijk summability bound \texttt{hZ}, the steepest-descent iterates satisfy
\[
\liminf_{k\to\infty}\ \left\|\nabla f(x_k)\right\| = 0,
\]
which for a non-negative sequence is exactly
\[
\forall\,\varepsilon > 0,\ \exists\,k,\ \texttt{vnorm}\,(\nabla f(x_k)) < \varepsilon.
\]
The proof rests on two auxiliary facts about the steepest-descent direction $d=-\nabla f(x)$. The lemma \texttt{cos\_theta\_neg\_grad\_eq\_one} shows that, whenever $\nabla f(x)\neq 0$, the angle cosine equals $1$; its squared form \texttt{cos\_theta\_neg\_grad\_sq\_eq\_one} gives $\cos^2\theta = 1$. The general principle is isolated as \texttt{gradient\_norm\_inf\_zero}: if the non-negative Zoutendijk terms $\cos^2\theta\,\|\nabla f\|^2$ have a uniform bound on every prefix sum, then no positive lower bound on $\|\nabla f\|$ can persist, so the gradient norm has infimum $0$. The positive statement is formalized as follows:
\begin{lstlisting}[
	basicstyle=\ttfamily\scriptsize,
	backgroundcolor=\color{gray!5},
	frame=single,
	mathescape=true,
	escapeinside={(*@}{@*)}
	]
(*@\textcolor{red}{lemma}@*) (*@\textcolor{blue}{cos\_theta\_neg\_grad\_eq\_one}@*)
  (x : Vec n) (hgx : vnorm (grad x) $\ne$ (*@\textcolor{green!50!black}{0}@*)) :
  cos_theta grad x (-grad x) = 1

(*@\textcolor{red}{lemma}@*) (*@\textcolor{blue}{cos\_theta\_neg\_grad\_sq\_eq\_one}@*)
  (x : Vec n) (hgx : vnorm (grad x) $\ne$ (*@\textcolor{green!50!black}{0}@*)) :
  (cos_theta grad x (-grad x))^2 = 1

(*@\textcolor{red}{lemma}@*) (*@\textcolor{blue}{gradient\_norm\_inf\_zero}@*)
  (xs : (*@$\mathbb{N}$@*) $\to$ Vec n) (L : (*@$\mathbb{R}$@*))
  (hZ : $\forall$ k,
    $\sum$ i $\in$ Finset.range k,
      cos_theta grad (xs i) (-grad (xs i)) ^ 2 *
      vnorm (grad (xs i)) ^ 2 $\le$ L) :
  $\forall$ $\varepsilon$ > (*@\textcolor{green!50!black}{0}@*), $\exists$ k, vnorm (grad (xs k)) < $\varepsilon$
\end{lstlisting}

\begin{lstlisting}[
	basicstyle=\ttfamily\scriptsize,
	backgroundcolor=\color{gray!5},
	frame=single,
	mathescape=true,
	escapeinside={(*@}{@*)}
	]
(*@\textcolor{red}{theorem}@*) (*@\textcolor{blue}{globalConvergence}@*)
  (c1 $\gamma$ $\alpha_0$ L : (*@$\mathbb{R}$@*))
  (hc1 : (*@\textcolor{green!50!black}{0}@*) < c1 $\land$ c1 < 1)
  (h$\gamma$ : (*@\textcolor{green!50!black}{0}@*) < $\gamma$ $\land$ $\gamma$ < 1)
  (h$\alpha_0$ : $\alpha_0$ > (*@\textcolor{green!50!black}{0}@*))
  (hdesc_all : $\forall$ x, isDescentDirection grad x (-grad x))
  (hfo : FirstOrderApprox f grad)
  (hZ : $\forall$ k,
    $\sum$ i $\in$ Finset.range k,
      cos_theta grad
        (x_seq f grad x0 c1 $\gamma$ $\alpha_0$
          hc1 h$\gamma$ h$\alpha_0$ hdesc_all hfo i)
        (-grad (x_seq f grad x0 c1 $\gamma$ $\alpha_0$
          hc1 h$\gamma$ h$\alpha_0$ hdesc_all hfo i)) ^ 2 *
      vnorm (grad (x_seq f grad x0 c1 $\gamma$ $\alpha_0$
        hc1 h$\gamma$ h$\alpha_0$ hdesc_all hfo i)) ^ 2 $\le$ L) :
  $\forall$ $\varepsilon$ > (*@\textcolor{green!50!black}{0}@*), $\exists$ k,
    vnorm (grad (x_seq f grad x0 c1 $\gamma$ $\alpha_0$
      hc1 h$\gamma$ h$\alpha_0$ hdesc_all hfo k)) < $\varepsilon$
\end{lstlisting}

\section*{Impact Statement}

This paper initiates a research direction on optimization formalization by providing a Lean 4 formalization of line search methods. Through strict machine-checked logic, we verify step-size update criteria, classical line search conditions, and key convergence results, thereby reducing ambiguity in textbook-style definitions and proofs. The work strengthens machine-verified foundation of nonlinear optimization and complements existing formalization of first-order methods, block-structured algorithms, and optimality conditions. Our work supports the transition of optimization theory from informal or empirical validation toward formalized and machine-verifiable reasoning. By establishing reusable formal components for line search, we lay a foundation for verifying more complex optimization algorithms and has significant value for both theoretical innovation and reliable engineering applications.

All code is available at \url{https://github.com/AquilaCheung/LEAN-formalization-lineSearchMethods/blob/main/zyy.lean}.


\begin{thebibliography}{18}
	\providecommand{\natexlab}[1]{#1}
	\providecommand{\url}[1]{\texttt{#1}}
	\expandafter\ifx\csname urlstyle\endcsname\relax
	\providecommand{\doi}[1]{doi: #1}\else
	\providecommand{\doi}{doi: \begingroup \urlstyle{rm}\Url}\fi
	
	\bibitem[Armijo(1966)]{armijo1966}
	L.~Armijo.
	\newblock Minimization of functions having lipschitz continuous first partial
	derivatives.
	\newblock \emph{Pacific Journal of Mathematics}, 16\penalty0 (1):\penalty0
	1--3, 1966.
	\newblock \doi{10.2140/pjm.1966.16.1}.
	
	\bibitem[{BICMR}()]{bicmr2026}
	{BICMR}.
	\newblock Optimization formalization platform.
	\newblock URL \url{http://faculty.bicmr.pku.edu.cn/~wenzw/formal/docs/#/}.
	
	\bibitem[Curry(1944)]{curry1944method}
	H.B. Curry.
	\newblock The method of steepest descent for nonlinear minimization problems.
	\newblock \emph{Quarterly of Applied Mathematics}, 2:\penalty0 258--263, 1944.
	
	\bibitem[Davidon(1959)]{davidon1959variable}
	W.C. Davidon.
	\newblock Variable metric method for minimization.
	\newblock Technical Report ANL-5990, Argonne National Lab, 1959.
	
	\bibitem[de~Moura and Ullrich(2021)]{demoura2021lean4}
	Leonardo de~Moura and Sebastian Ullrich.
	\newblock The lean 4 theorem prover and programming language.
	\newblock In \emph{Automated Deduction -- CADE 28: 28th International
		Conference on Automated Deduction, Virtual Event, July 12--15, 2021,
		Proceedings}, pages 625--635, Berlin, Heidelberg, 2021. Springer.
	\newblock ISBN 978-3-030-79875-8.
	\newblock \doi{10.1007/978-3-030-79876-5\_37}.
	\newblock URL \url{https://doi.org/10.1007/978-3-030-79876-5\_37}.
	
	\bibitem[Fletcher and Powell(1963)]{fletcher1963rapidly}
	R.~Fletcher and M.J.D. Powell.
	\newblock A rapidly convergent descent method for minimization.
	\newblock \emph{The Computer Journal}, 6:\penalty0 163--168, 1963.
	
	\bibitem[Goldstein(1965)]{goldstein1965}
	A.~A. Goldstein.
	\newblock On steepest descent.
	\newblock \emph{SIAM Journal on Control}, 3\penalty0 (1):\penalty0 147--151,
	1965.
	\newblock \doi{10.1137/0303013}.
	
	\bibitem[Grippo et~al.(1986)Grippo, Lampariello, and
	Lucidi]{grippo1986nonmonotone}
	Luigi Grippo, Francesco Lampariello, and Stefano Lucidi.
	\newblock A nonmonotone line search technique for {N}ewton's method.
	\newblock \emph{SIAM Journal on Numerical Analysis}, 23\penalty0 (4):\penalty0
	707--716, 1986.
	\newblock \doi{10.1137/0723046}.
	
	\bibitem[Kiefer(1953)]{kiefer1953sequential}
	Jack Kiefer.
	\newblock Sequential minimax search for a maximum.
	\newblock \emph{Proceedings of the American Mathematical Society}, 4\penalty0
	(3):\penalty0 502--506, 1953.
	
	\bibitem[Li et~al.(2024)Li, Wang, He, Wu, Xu, and Wen]{li2024formalization}
	Chenyi Li, Ziyu Wang, Wanyi He, Yuxuan Wu, Shengyang Xu, and Zaiwen Wen.
	\newblock Formalization of complexity analysis of the first-order algorithms
	for convex optimization.
	\newblock \emph{arXiv preprint arXiv:2403.11437}, 2024.
	
	\bibitem[Li et~al.(2025{\natexlab{a}})Li, Wang, Bai, Duan, Hao, Gao, and
	Wen]{li2025block}
	Chenyi Li, Zichen Wang, Yifan Bai, Yunxi Duan, Pengfei Hao, Yuqing Gao, and
	Zaiwen Wen.
	\newblock Formalization of algorithms for optimization with block structures.
	\newblock \emph{arXiv preprint arXiv:2503.18806}, 2025{\natexlab{a}}.
	
	\bibitem[Li et~al.(2025{\natexlab{b}})Li, Xu, Sun, Zhou, and Wen]{li2025kkt}
	Chenyi Li, Shengyang Xu, Chumin Sun, Li~Zhou, and Zaiwen Wen.
	\newblock Formalization of optimality conditions for smooth constrained
	optimization problems.
	\newblock \emph{arXiv preprint arXiv:2503.18821}, 2025{\natexlab{b}}.
	
	\bibitem[mathlib Community(2020)]{mathlib2020}
	The mathlib Community.
	\newblock The lean mathematical library.
	\newblock In \emph{Proceedings of the 9th ACM SIGPLAN International Conference
		on Certified Programs and Proofs (CPP)}. ACM, 2020.
	\newblock \doi{10.1145/3372885.3373824}.
	
	\bibitem[Wolfe(1969)]{wolfe1969}
	P.~Wolfe.
	\newblock Convergence conditions for ascent methods.
	\newblock \emph{SIAM Review}, 11\penalty0 (2):\penalty0 226--235, 1969.
	\newblock \doi{10.1137/1011036}.
	
	\bibitem[Wolfe(1971)]{wolfe1971}
	P.~Wolfe.
	\newblock Convergence conditions for ascent methods {II}: Some corrections.
	\newblock \emph{SIAM Review}, 13\penalty0 (2):\penalty0 185--188, 1971.
	\newblock \doi{10.1137/1013035}.
	
	\bibitem[Yang(2012)]{yang2012}
	Y.~Yang.
	\newblock A robust {BFGS} algorithm for unconstrained nonlinear optimization
	problems, 2012.
	\newblock URL \url{https://arxiv.org/pdf/1212.5929}.
	\newblock arXiv preprint arXiv:1212.5929.
	
	\bibitem[Zhang and Hager(2004)]{zhang2004nonmonotone}
	Houchun Zhang and William~W. Hager.
	\newblock A nonmonotone line search technique and its application to
	unconstrained optimization.
	\newblock \emph{SIAM Journal on Optimization}, 14\penalty0 (4):\penalty0
	1043--1056, 2004.
	\newblock \doi{10.1137/S1052623403428208}.
	
	\bibitem[Zoutendijk(1960)]{zoutendijk1960}
	G.~Zoutendijk.
	\newblock \emph{Methods of feasible directions: A study in linear and
		non-linear programming}.
	\newblock Elsevier, 1960.
	
\end{thebibliography}
	
	
\end{document}